\documentclass{amsart}
\usepackage{cases}
\usepackage[dvips]{graphicx}
\usepackage{mathrsfs}
\usepackage{float}
\usepackage{verbatim}

\newcommand{\Z}{{\mathbb{Z}}}

\newcommand{\R}{{\mathbb{R}}}
\newcommand{\C}{{\mathbb{C}}}

\theoremstyle{definition}
  \newtheorem{Thm}{Theorem}[section]
  \newtheorem{Lem}[Thm]{Lemma}
  \newtheorem{Def}[Thm]{Definition}
  \newtheorem{Ex}[Thm]{Example}
  
  \newtheorem{Prop}[Thm]{Proposition}

\makeatletter
\theoremstyle{remark}
\newenvironment{pro}{\begin{proof}}{\end{proof}}

\thanks{This paper contains some results of my master's thesis contributed to Tokyo Institute of Technology.}
\date{\today}
\begin{document}
\title{On the Conway potential function introduced by Kauffman}
\author{Masashi Sato}
\address{
Department of Mathematics,
Tokyo Institute of Technology,
Oh-okayama, Meguro, Tokyo 152-8551, Japan
}
\email{09m00114@math.titech.ac.jp}
\maketitle
\begin{abstract}
 We show two results about the Conway potential function which is known as the normalized multivariable Alexander polynomial. 
 We first show that the Conway potential function introduced by Kauffman in ``Formal Knot Theory'' is indeed a link invariant. Next we show that Kauffman's potential function equals Hartley's potential function. We will prove it by using Murakami's axioms for the multivariable Alexander polynomial. 
\end{abstract}
\section{Introduction}
 In \cite{MR0258014} J.H.~Conway introduced the Conway potential function by some axioms. But his axioms are not sufficient to determine his potential function. L.H.~Kauffman \cite{MR592573} showed how to define the one variable reduced potential function in terms of a Seifert matrix. R.~Hartley \cite{MR727708} gave a precise definition of the multivariable potential function. At the same time Kauffman \cite{MR712133} introduced another definition of the multivariable potential function without using Seifert matrices. However he did not show that his definition gave a link invariant. In this paper we show it in Section~3. 

On the other hand J.~Murakami \cite{MR1177442,MR1197048} gave axioms of the potential function which are sufficient for the definition. In Section~4 we show that Kauffman's potential function equals Hartley's potential function by J.~Murakami's axioms. In Section~2 we confirm some necessary definitions and theorems in \cite{MR712133}.
\section{preparation}
In this section we introduce some definitions and theorems in \cite{MR712133} to define Kauffman's potential function. 

In this paper we regard a link diagram as a regular projection of a link with over and under informations at the vertices. On the other hand a link projection is without such information. 
\begin{Def}[Kauffman state \cite{MR712133}]
Let $U$ be a link projection on $\R ^2$ or $S ^2$. A pair of a vertex of $U$ and a local region around the vertex is represented as a marker as in Figure~\ref{fig:mark}. A local region means intersection of a region and a interior of the dotted circle. 
\begin{figure}[h]
  \begin{center}
  \includegraphics[scale=1.25]{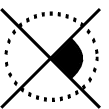}
  \caption{marker}
  \label{fig:mark}
  \end{center}
\end{figure}

If $U$ is connected, the number of regions of $U$ is two more than the number of the vertices since the Euler characteristic of $S ^2$ is two. 
Then we call a set of markers a {\em Kauffman state} or a {\em state} if it satisfies the following conditions.
 \begin{enumerate}
   \item Two adjacent regions have no markers. We put a star ($\ast$) in each of these regions. 
   \item Every vertex has only one marker.
   \item Every region without a star ($\ast$) contains only one marker.  
 \end{enumerate}
\end{Def}
\begin{Ex}Figure~\ref{fig:Kauffman state} shows an example of a knot projection and a state. 
\begin{figure}[h]
 \begin{center}
 \includegraphics[scale=0.8]{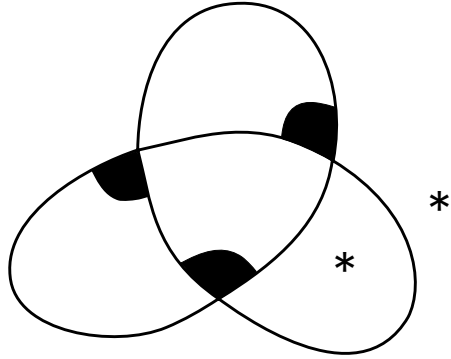}
 \caption{Kauffman state}
 \label{fig:Kauffman state}
 \end{center}
\end{figure}
\end{Ex}
\begin{Def}[transposition \cite{MR712133}]
The operation as shown in  Figure~\ref{fig:tp} is called a {\em transposition}. In particular the operation from left to right is called a {\em clockwise transposition}, and that from right to left is called a {\em counter-clockwise transposition}. The dotted circle contains a part of a projection and (possibly empty) markers. Moreover a state is said to be a {\em clocked state} if it admits only clockwise transpositions and a {\em counter-clocked state} if it admits only counter-clockwise transpositions. 
\begin{figure}[h]
 \begin{center}
 \includegraphics{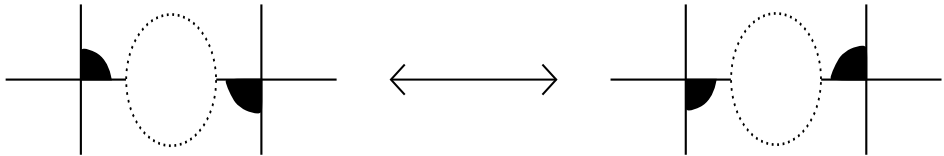}
 \caption{transposition}
 \label{fig:tp}
 \end{center}
\end{figure}
\end{Def}
\begin{Thm}[Clock Theorem \cite{MR712133}]
Let $U$ be a connected link projection and $\mathscr{S}$ be the set of states of $U$ for a given choice of adjacent stars. 

Then the set $\mathscr{S}$ has the following properties. 
\begin{enumerate}
 \item It has a unique clocked state and a unique counter-clocked state. 
 \item Any state in $\mathscr{S}$ can be reached from the clocked (counter-clocked respectively) state by a series of clockwise (counter-clockwise respectively) transpositions. 
\end{enumerate}
\end{Thm}
\begin{Def}[state polynomial \cite{MR712133}]
Let $U$ be an  oriented link projection with $n$ vertices. Let $\mathscr{S}$ be the set of states of $U$ with fixed adjacent stars. We label the vertices as $v_1,v_2,\dots ,v_n$ and put variables $I_1,O_1,U_1,D_1,\dots ,I_n,O_n,U_n,D_n$ in local regions around the vertices as indicated in Figure~\ref{fig:v_v}. Then we define $\langle U|S\rangle$ for $S \in \mathscr{S}$ by the formula $\langle U|S\rangle = V_1(S)V_2(S)\dots V_n(S)$ where $V_k(S)$ is the variable touched by the marker defined by $S$ at the vertex $v_k$. 
\begin{figure}[h]
 \begin{center}
 \includegraphics{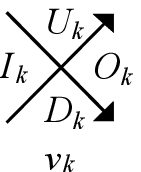}
 \caption{}
 \label{fig:v_v}
 \end{center}
\end{figure}
\\ We define $\sigma (S)$ to be $\langle U|S \rangle$ replacing all the $I_k$ with~$-1$ and the $O_k$, $U_k$ and $D_k$ with~$1$. We regard $\langle U|S\rangle $ as an element of $\C [I_1,O_1,U_1,D_1,\dots ,I_n,O_n,U_n,D_n]$. If $U$ is connected, then the {\em state polynomial} for $U$ is defined by the formula: 
\begin{equation*}
\langle U| \mathscr{S} \rangle =\sum_{S\in \mathscr{S}} \sigma (S)\langle U|S \rangle \in \C [I_1,O_1,U_1,D_1,\dots ,I_n,O_n,U_n,D_n]. 
\end{equation*}
If $U$ is non-connected, we define $\langle U| \mathscr{S} \rangle = 0$. 

Moreover let $\langle U| \mathscr{S} \rangle ' \in \C [B,W]$ be obtained from $\langle U| \mathscr{S} \rangle$ replacing all the $I_k$, $O_k$, $U_k$ and $D_k$ with $1$, $B$ or $W$ as indicated in Figure~\ref{fig:bw0} .
\begin{figure}[h]
 \begin{center}
 \includegraphics{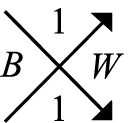}
 \caption{}
 \label{fig:bw0}
 \end{center}
\end{figure}
\end{Def}
\begin{Def}[Alexander matrix \cite{MR712133}]
Let $U$ be a connected link projection with $n$ vertices. We label the vertices and the regions as $v_1,v_2,\dots ,v_n$ and $r_1,r_2,\dots ,r_{n+2}$ respectively. We put variables as indicated in Figure~\ref{fig:v_v}. Then we define the {\em Alexander matrix} $A(U)$ as an $n \times (n+2)$ matrix with $(i,j)$ entry $A(U)_{ij}$, where $A(U)_{ij}$ is the sum of the variables around the vertex $v_i$ in the region $r_j$. Let $A(U)(i_1,i_2)$ be the $n \times n$ matrix obtained from $A(U)$ by deleting the $i_1$th and $i_2$th columns. It is called the {\em reduced Alexander matrix}. Moreover let $A'(U)$ be obtained from $A(U)$ by replacing all the $I_k$, $O_k$, $U_k$ and $D_k$ with $1$, $B$ or $W$ as indicated in Figure~\ref{fig:bw0} .
\end{Def}
\begin{Thm}[\cite{MR712133}]\label{thm:sp_am}
Let $U$ be a connected link projection with labeled vertices and regions as $v_1,v_2,\dots ,v_n$ and $r_1,r_2,\dots ,r_{n+2}$ respectively. Let $\mathscr{S}(i_1,i_2)$ be the set of states of $U$ with fixed adjacent stars in the regions $r_{i_1}$ and $r_{i_2}$. Then state polynomial and the Alexander matrix satisfy the following formula:
\begin{equation*}
\langle U| \mathscr{S}(i_1,i_2) \rangle \doteq \det A(U)(i_1, i_2) .
\end{equation*}
The symbol $\doteq$ means equality up to sign.  
\end{Thm}
\begin{Def}[multiple Alexander index \cite{MR712133}]
Let $L=L_1 \cup L_2 \cup \dots \cup L_N$ be an oriented link diagram on $\R ^2$. Let $U$ be the projection of $L$. Then each region of $U$ is assigned an element $p=(p_1,p_2,\dots,p_N)$ of $\underbrace{\Z \times \Z \times \dots  \times \Z}_{N}$ as follows. The unbounded region is assigned index $(0,0,\dots ,0)$. The $K$th component, $p_K$, increases if one crosses the $K$th strand from left to right and decrease if one crosses the $K$th strand from right to left. 
Then each element $p=(p_1,p_2,\dots,p_N)$ is called the {\em multiple Alexander index} of $U$. Moreover the sum of the components of the multiple Alexander index in a region is called the ({\em non-multiple}) {\em Alexander index}. 
\begin{figure}[h]
 \begin{center}
 \includegraphics{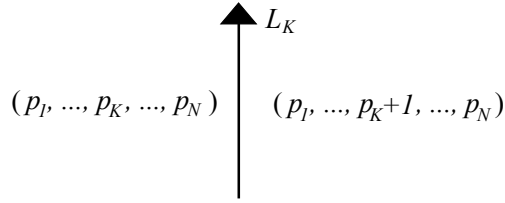}
 \caption{multiple Alexander index}
 \label{fig:mai}
 \end{center}
\end{figure} 
\end{Def}
The following lemma is proved in \cite{MR712133}. Since we use the technique in the proof later, we give a proof following Kauffman. 
\begin{Lem}[\cite{MR712133}]
Let $U$ be an oriented link projection. Let $\mathscr{S}_1$ and $\mathscr{S}_2$ be two sets of states with adjacent stars. Then the state polynomials satisfy the formula $\langle U| \mathscr{S}_1 \rangle ' = \langle U| \mathscr{S}_2 \rangle '$. 
\end{Lem}
\begin{pro}
If $U$ is not connected, $\langle U| \mathscr{S}_1 \rangle ' = 0 = \langle U| \mathscr{S}_2 \rangle '$. So we assume that $U$ is connected. Now we give the non-multiple Alexander index to $U$, and label the vertices and regions as $v_1,v_2,\dots ,v_n$ and $r_1,r_2,\dots ,r_{n+2}$ respectively. Let $m(j)$ be the Alexander index in the region $r_j$. 

First let $\alpha $ and $\alpha ^{-1}$ be the roots of $X^2 + (B+W)X +1 = 0$, where $X$ is a variable and $B$ and $W$ are regarded as constants. Let $\beta (x, y) = \alpha ^{ m(x) - m(y) } - \alpha ^{-m(x) + m(y)} $ and $\dot{a}_j$ be the $j$th column of the Alexander matrix $A'(U)$. Then we have 
\begin{equation}
 \sum_{j=1}^{n+2} \beta (j,k) \dot{a}_j = \bf{0} \label{eq:b1}
\end{equation}
for a fixed integer $k$ from the relation between the Alexander index and the Alexander matrix (See~p.58~of~\cite{MR712133}). The symbol $\bf{0}$ means a zero-vector. Therefore we have
\begin{equation}
 \begin{split}
  &\beta (l,j) \det A'(U)(j,k) \\
                             &= \det \begin{pmatrix} \dot{a}_1 & \dots & \hat{\dot{a}}_j & \dots & \hat{\dot{a}}_k & \dots & \beta (l,j) \dot{a}_l & \dots & \dot{a}_{n+2} \end{pmatrix} \\
                             &= \det \begin{pmatrix} \dot{a}_1 & \dots & \hat{\dot{a}}_j & \dots & \hat{\dot{a}}_k & \dots & -\sum_{j' \neq l} \beta (j',j) \dot{a}_{j'} & \dots & \dot{a}_{n+2} \end{pmatrix} \\
                             &= \det \begin{pmatrix} \dot{a}_1 & \dots & \hat{\dot{a}}_j & \dots & \hat{\dot{a}}_k & \dots & -\beta (j,j) \dot{a}_{j} - \beta (k,j) \dot{a}_k & \dots & \dot{a}_{n+2} \end{pmatrix} \\
                             &= - \beta (k,j) \det \begin{pmatrix} \dot{a}_1 & \dots & \hat{\dot{a}}_j & \dots & \hat{\dot{a}}_k & \dots & \dot{a}_k & \dots & \dot{a}_{n+2} \end{pmatrix} \\
                            & \doteq  \beta (k,j) \det A'(U)(j,l)
 \end{split}
\label{eq:b2}
\end{equation}
for any fixed $j$, $k$ and $l$ so that $j \neq k$ and $j \neq l$. Similarly we have 
\begin{equation*}
\beta (h,l) \det  A'(U)(l,j) \doteq \beta (j,l) \det A'(U)(l,h)
\end{equation*}
for any fixed $h$ so that $l \neq h$. Therefore we have the following formula: 
\begin{equation}
 \frac{\beta (l,j)}{\beta (k,j)} \det A'(U)(j,k) \doteq \frac{\beta (j,l)}{\beta (h,l)} \det A'(U)(l,h).
\end{equation}
Now we have $\beta (l,j) = -\beta (j,l)$ from the definition of $\beta$. So we have 
\begin{equation*}
\frac{\det A'(U)(j,k)}{\beta (k,j)} \doteq \frac{\det A'(U)(l,h)}{\beta (h,l)}.
\end{equation*}
Now we assume that the regions $r_{j}$ and $r_{k}$ are adjacent and that $r_{l}$ and $r_{h}$ are adjacent. Then we have the formulas $\det A'(U)(j, k) \doteq \langle U| \mathscr{S}(j,k) \rangle '$ and $\det A'(U)(l, h) \doteq \langle U| \mathscr{S}(l,h) \rangle '$ from Theorem \ref{thm:sp_am}. Moreover we have $\beta (k,j) \doteq \alpha - \alpha ^{-1} \doteq \beta (h,l)$. 
Therefore we have $\langle U| \mathscr{S}(j,k) \rangle ' \doteq \langle U| \mathscr{S}(l,h) \rangle '$. On the other hand we have $\langle U| \mathscr{S} \rangle ' = \sum_{S \in \mathscr{S}}(-1)^{b(S)}B^{b(S)}W^{w(S)}$ where $b(S)$ denotes the number of $B$ in $S$, and $w(S)$ denotes the number of $W$ in $S$ from the definition. Therefore we have $\langle U| \mathscr{S}(j,k) \rangle ' = \langle U| \mathscr{S}(l,h) \rangle '$. 
\end{pro}
Now we prepare the following propositions to prove the invariance of Kaffman's potential function later. 
\begin{Prop}\label{pro:sp_am}
Let $U$ be an oriented and connected link projection with labeled vertices and regions as $v_1,v_2,\dots ,v_n$ and $r_1,r_2,\dots ,r_{n+2}$ respectively so that the regions $r_{n+1}$ and $r_{n+2}$ are adjacent as indicated in Figure~\ref{fig:as0} and that the Alexander matrix $A(U)$ satisfies $\langle U| \mathscr{S} (n+1,n+2) \rangle = \det A(U)(n+1, n+2)$. Then the Alexander matrix and the state polynomial satisfy the following formulas: 
\begin{subnumcases}
{\det A(U)(i,j)=}
(-1)^{i+j+1} \langle U| \mathscr{S} (i,j) \rangle \text{\ \ \ \ \ (Figure~\ref{fig:as1})} \label{eq:as_a}\\
(-1)^{i+j}\langle U| \mathscr{S} (i,j) \rangle \text{\ \ \ \ \ \ \ \ (Figure~\ref{fig:as2})} \label{eq:as_b},
\end{subnumcases}
where $j > i$.  
\begin{figure}[h]
 \begin{minipage}{0.3\hsize}
  \begin{center}
   \includegraphics{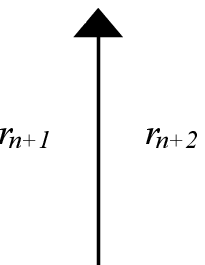}
  \end{center}
  \caption{}
  \label{fig:as0}
 \end{minipage}
 \begin{minipage}{0.3\hsize}
  \begin{center}
   \includegraphics{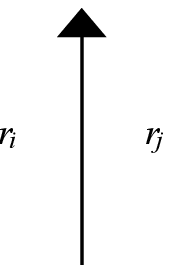}
  \end{center}
  \caption{}
  \label{fig:as1}
 \end{minipage}
 \begin{minipage}{0.3\hsize}
  \begin{center}
   \includegraphics{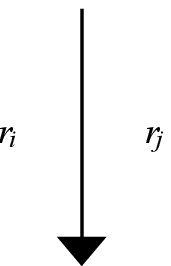}
  \end{center}
  \caption{}
  \label{fig:as2}
 \end{minipage}
\end{figure}
\end{Prop}
In order to prove this proposition we prepare the following lemma. 
\begin{Lem}\label{lem:am_am}
We have $\det A(U)(i,j) = (-1)^a \langle U| \mathscr{S} (i,j) \rangle $ if $\det A'(U)(i,j) = (-1)^a \langle U| \mathscr{S} (i,j) \rangle '$ for some integer $a$. 
\end{Lem}
\begin{pro}
We can prove this since $\langle U| \mathscr{S} \rangle ' = \sum_{S \in \mathscr{S}}(-1)^{b(S)}B^{b(S)}W^{w(S)} \neq 0$ for a non-empty set $\mathscr{S}$. 
\end{pro}
\begin{pro}[Proof of Proposition \ref{pro:sp_am}]
We consider $A'(U)$ instead of $A(U)$ in this proof from Lemma~\ref{lem:am_am}. We replace $j$ with $n+2$, $k$ with $n+1$, $l$ with $i$ and $h$ with $j$ in the proof of Theorem~\ref{thm:sp_am} and we discuss the sign. Then from~\eqref{eq:b2} we have 
\begin{equation}
 \begin{split}
  &\beta (i,n+2) \det A'(U)(n+2,n+1) \\
                             &= \det \begin{pmatrix} \dot{a}_1 & \dots & \beta (i,n+2) \dot{a}_i & \dots &  \hat{\dot{a}}_{n+1} & \hat{\dot{a}}_{n+2} \end{pmatrix} \\
                             &= -\beta (n+1,n+2) \det \begin{pmatrix} \dot{a}_1 & \dots & \dot{a}_{n+1} & \dots &  \hat{\dot{a}}_{n+1} & \hat{\dot{a}}_{n+2} \end{pmatrix} \\
                             &= (-1)^{n-i+1} \beta (n+1,n+2) \det A'(U)(i,n+2),
 \end{split}
\label{eq:b3}
\end{equation} 
where $i \leq n+1$. Similarly we have 
\begin{equation}
 \begin{split}
  &\beta (j,i) \det A'(U)(i,n+2) \\
                             &= \det \begin{pmatrix} \dot{a}_1 & \dots & \hat{\dot{a}}_i & \dots & \beta (j,i) \dot{a}_j & \dots & \hat{\dot{a}}_{n+2} \end{pmatrix} \\
                             &= - \beta (n+2,i) \det \begin{pmatrix} \dot{a}_1 & \dots & \hat{\dot{a}}_i & \dots & \dot{a}_{n+1} & \dots & \hat{\dot{a}}_{n+2} \end{pmatrix} \\
                             &= (-1)^{n-j} \beta (n+2,i) \det A'(U)(i,j).
 \end{split}
\label{eq:b4}
\end{equation}
Therefore we have 
\begin{equation}
 \frac{\beta (i,n+2)}{\beta (n+1,n+2)} \det A'(U)(n+2,n+1) = (-1)^{i+j+1} \frac{\beta (n+2,i)}{\beta (j,i)} \det A'(U)(i,j). 
\label{eq:b5}
\end{equation}
On the other hand we have $\beta (i,n+2) = - \beta (n+2,i)$ and $\beta (n+1,n+2)= \alpha ^{-1} - \alpha $ from the definitions of $\beta$ and the Alexander index. If the regions $r_i$ and $r_j$ are as shown in Figure~\ref{fig:as1}, we have $\beta (j,i)= \alpha - \alpha ^{-1}$. So we have $\det A'(U)(n+2,n+1) = (-1)^{i+j+1} \det A'(U)(i,j)$ from \eqref{eq:b5}. Then we have 
\begin{equation*}
\begin{split}
\langle U| \mathscr{S} (i,j) \rangle ' &= \langle U| \mathscr{S} (n+1,n+2) \rangle ' \\
                                       &= \det A'(U)(n+2,n+1) \\
                                       &= (-1)^{i+j+1} \det A'(U)(i,j). 
\end{split}
\end{equation*}
Similarly if the regions $r_i$ and $r_j$ are as shown in Figure~\ref{fig:as2}, we can get \eqref{eq:as_b}. This proof is complete. 
\end{pro}
%
%
\section{Kauffman's potential function}
In this section we will define the potential function introduced by Kauffman and discuss its invariance. 
\begin{Def}[potential function \cite{MR712133}]
Let $L=L_1 \cup L_2 \cup \dots \cup L_N$ be an oriented link diagram. Let $U$ be the projection of $L$ with labeled vertices as $v_1,v_2,\dots ,v_n$ and the multiple Alexander index. Let $\mathscr{S}$ be the set of states of $U$ whose fixed stars share the $K$th strand with indices $(p_1, \dots , p_K , \dots ,p_N)$ and $(p_1, \dots , p_K +1, \dots ,p_N)$. First we put $N$~variables $X_1, X_2, \dots , X_N$ in local regions around the vertices as shown in Figure~\ref{fig:pf_1}. Then let $\langle L| \mathscr{S} \rangle$ be the polynomial obtained from $\langle U| \mathscr{S} \rangle$ by replacing all the $I_k$, $O_k$, $D_k$ and $U_k$ with $X_1, X_2, \dots , X_N$ as indicated in Figure~\ref{fig:pf_1}. Second let $|\mathscr{S}|$ be $X_1^{-2p_1} X_2^{-2p_2} \cdots X_N ^{-2p_N }X_K ^{-1}(X_K - X_K ^{-1})$. Third let $c(L_J)$ be the curvature of the sublink $L_J$, which counts how many times the sublink rotates counter-clockwise. 

Then we define $\square _L$ to be
\begin{equation*}
 \square _L = \frac{X_1^{c(L_1)} X_2^{c(L_2)} \cdots X_N ^{c(L_N) }}{\left | \mathscr{S} \right |} \langle L| \mathscr{S} \rangle 
\end{equation*}
and we call $\square _L$ {\em  Kauffman's potential function}. 
\begin{figure}[h]
 \begin{center}
 \includegraphics{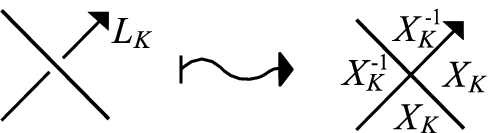}
 \caption{}
 \label{fig:pf_1}
 \end{center}
\end{figure}
\end{Def}
\begin{Thm}\label{pro:pf}
Kauffman's potential function is a link invariant. 
\end{Thm}
In order to prove this theorem we prove the following two lemmas. 
\begin{Lem}\label{lem:pf1}
The function $\langle L| \mathscr{S} \rangle / |\mathscr{S}|$ is independent of the choice of fixed stars.
\end{Lem}
\begin{pro}
Let $U$ be the projection of $L$ with labeled vertices and regions as $v_1,v_2,\dots ,v_n$ and $r_1,r_2,\dots ,r_{n+2}$ respectively so that they satisfy the condition of Proposition~\ref{pro:sp_am}. Let $A(L) = (\dot{a}_1, \dots , \dot{a}_{n+2})$ be the Alexander matrix corresponding to Figure~\ref{fig:pf_1}. Let $m_K(j)$ be the $K$th component of the multiple Alexander index in the region $r_j$. 

First we show the following two formulas: 
\begin{subequations}
\begin{align}
&\sum_{j=1}^{n+2} (-1)^{ \{ m_1(j)+m_2(j)+ \dots +m_N (j) \} }\dot{a}_j = \bf{0} \label{eq:pf_a},  \\
&\sum_{j=1}^{n+2} (-1)^{ \{ m_1(j)+m_2(j)+ \dots +m_N (j) \} }  X_1^{-2m_1(j)}X_2^{-2m_2(j)} \cdots X_N ^{-2m_N (j)} \dot{a}_j = \bf{0} \label{eq:pf_b}. 
\end{align}
\end{subequations}
Let $a_{ij}$ be the $(i,j)$ entry of $A(U)$. For any vertex $v_i$ with a positive crossing we have the following two formulas from Figure~\ref{fig:pf_2}, where $(\delta _i(1),\dots , \delta _i(N))$ is the multiple Alexander index in the region touching the vertex $v_i$ such that the region is above the vertex: 
\begin{equation*}
\begin{split}
&\sum_{j=1}^{n+2} (-1)^{ \{ m_1(j)+m_2(j)+ \dots +m_N (j) \} }a_{ij} \\
&= (-1)^{\{ \delta _{i}(1)+\delta _{i}(2)+ \dots +\delta _{i}(N)  \}} \{ (-1)^2X_K + (-1)(X_K^{-1} + X_K) + X_K ^{-1}\} \\ 
&= 0, \\
&\text{and} \\
&\sum_{j=1}^{n+2} (-1)^{ \{ m_1(j)+m_2(j)+ \dots +m_N (j) \} }  X_1^{-2m_1(j)}X_2^{-2m_2(j)} \cdots X_N ^{-2m_N (j)} a_{ij} \\ 
&= (-1)^{\sum_{I =1}^{N }\delta _{i}(I) } \prod_{I =1}^{N }X_I ^{-2\delta _{i}(I) } \{ (-1)^2X_K ^{-2}X_J ^{-2}X_K + (-X_J ^{-2})X_K ^ {-1} + (-X_K ^{-2})X_K  + X_K ^{-1}\}  \\ 
&= 0 .
\end{split}
\end{equation*}
\begin{figure}[H]
 \begin{center}
 \includegraphics[scale=0.7]{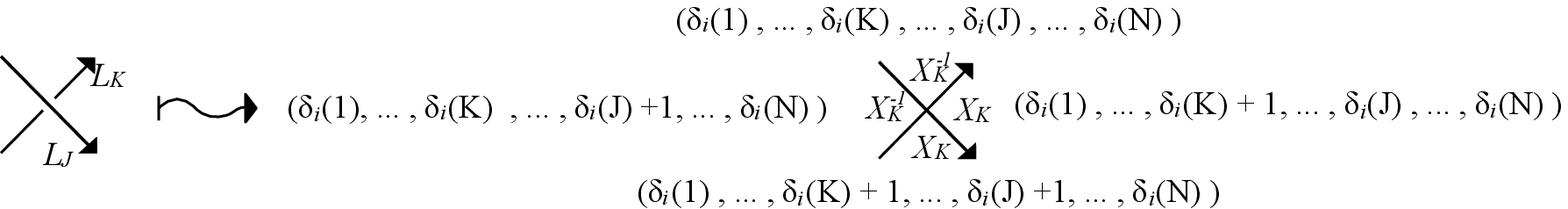}
 \caption{}
 \label{fig:pf_2}
 \end{center}
\end{figure}
On the other hand for any vertex $v_i$ with a negative crossing we have the following two formulas from Figure~\ref{fig:pf_3}: 
\begin{equation*}
\begin{split}
&\sum_{j=1}^{n+2} (-1)^{ \{ m_1(j)+m_2(j)+ \dots +m_N (j) \} }a_{ij} \\
&= (-1)^{\{ \delta _{i}(1)+\delta _{i}(2)+ \dots +\delta _{i}(N)  \}} \{ (-1)^2X_J + (-1)(X_J + X_J^{-1} ) + X_J^{-1} \} \\ 
&= 0 , \\
&\text{and} \\
&\sum_{j=1}^{n+2} (-1)^{ \sum_{I =1}^{N } m_I (j) }\prod_{I =1}^{N }X_I  ^{-2m_I (j) } a_{ij} \\ 
&= (-1)^{\sum_{I =1}^{N }\delta _{i}(I) } \prod_{I =1}^{N }X_I ^{-2\delta _{i}(I) } \{ (-1)^2X_K ^{-2}X_J ^{-2}X_J + (-X_J ^{-2})X_J+ (-X_J ^{-2})X_J^{-1} + X_J ^{-1}\} \\ 
&= 0.
\end{split}
\end{equation*}
\begin{figure}[H]
 \begin{center}
 \includegraphics[scale=0.7]{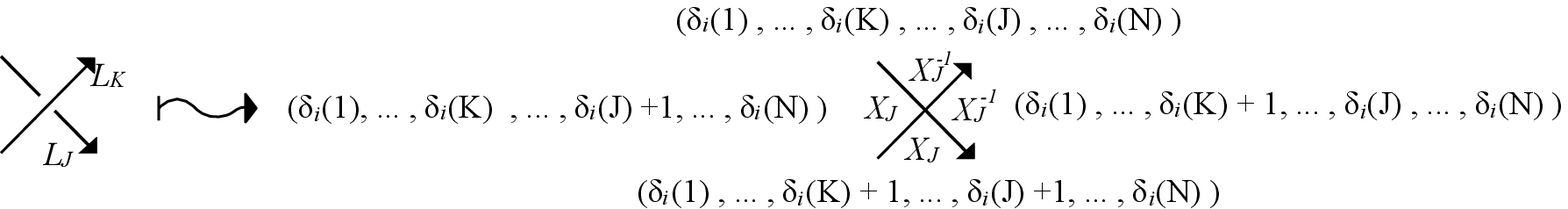}
 \caption{}
 \label{fig:pf_3}
 \end{center}
\end{figure}
Therefore we get \eqref{eq:pf_a} and \eqref{eq:pf_b}. From \eqref{eq:pf_a} and \eqref{eq:pf_b} we have the following formulas: 
\begin{equation*}
\sum_{j=1}^{n+2} (-1)^{ \sum_{I =1}^{N } \{ m_I (j)-m_I (k)\} }\dot{a}_j = \bf{0} ,
\end{equation*}
and
\begin{equation*}
\sum_{j=1}^{n+2} (-1)^{ \sum_{I =1}^{N } \{ m_I (j)-m_I (k)\} } \prod_{I =1}^{N }X_I ^{-2 \{ m_I (j)-m_I (k)\} } \dot{a}_j = \bf{0} 
\end{equation*}
for any $k$. So letting $\beta '(x,y) = (-1)^{ \sum_{I =1}^{N } \{ m_I (x)-m_I (y)\}} \{ 1- \prod_{I =1}^{N }X_I ^{-2 \{ m_I (x)-m_I (y)\} } \}$, we have
\begin{equation}
 \sum_{j=1}^{n+2} \beta '(j,k) \dot{a}_j = \bf{0} \label{eq:pf_1}. 
\end{equation} 
Then we can get the following formula since \eqref{eq:pf_1} is equal to \eqref{eq:b1}: 
\begin{equation*}
 \frac{\beta '(i,n+2)}{\beta '(n+1,n+2)} \det A(L)(n+2,n+1) = (-1)^{i+j+1} \frac{\beta '(n+2,i)}{\beta '(j,i)} \det A(L)(i,j).
\end{equation*}
Moreover we have \begin{equation*}
\begin{split}
 \frac{\beta '(i,n+2)}{\beta '(n+2,i)} &= \frac{1- \prod_{I =1}^{N }X_I ^{-2 \{ m_I (i)-m_I (n+2)\} }}{1- \prod_{I =1}^{N }X_I ^{-2 \{ m_I (n+2)-m_I (i)\} }}  \\ 
&= -\frac{X_1^{-2m_1(i)}X_2^{-2m_2(i)} \cdots X_N ^{-2m_N (i)}}{X_1^{-2m_1(n+2)}X_2^{-2m_2(n+2)} \cdots X_N ^{-2m_N (n+2)}} \\
 &= -X_J ^{-2} \frac{X_1^{-2m_1(i)}X_2^{-2m_2(i)} \cdots X_N ^{-2m_N (i)}}{X_1^{-2m_1(n+1)}X_2^{-2m_2(n+1)} \cdots X_N ^{-2m_N (n+1)}},
\end{split}
\end{equation*}
where the regions $r_{n+1}$ and $r_{n+2}$ are as shown in Figure~\ref{fig:as0} with the strand labeled as $J$. So we have 
\begin{equation*}
\begin{split}
 &-X_J ^{-2} \frac{\det A(L)(n+2,n+1)}{\beta '(n+1,n+2)X_1^{-2m_1(n+1)}X_2^{-2m_2(n+1)} \cdots X_N ^{-2m_N (n+1)}} \\
=& (-1)^{i+j+1} \frac{\det A(L)(i,j)}{\beta '(j,i)X_1^{-2m_1(i)}X_2^{-2m_2(i)} \cdots X_N ^{-2m_N (i)}} .
\end{split}
\end{equation*}
Then we have $\beta ' (n+1,n+2) = -(1-X_J ^{-2})$ and if the regions $r_i$ and $r_j$ are as shown in Figure~\ref{fig:as1} with the strand labeled as $M$, we have $\beta ' (j,i) = -(1-X_M ^{-2})$. Therefore we have 
\begin{align*}
 &\frac{\det A(L)(n+2,n+1)}{X_1^{-2m_1(n+1)}X_2^{-2m_2(n+1)} \cdots X_N ^{-2m_N (n+1)}X_J ^{-1}(X_J - X_J ^{-1})}  \\
=& (-1)^{i+j+1} \frac{\det A(L)(i,j))}{X_1^{-2m_1(i)}X_2^{-2m_2(i)} \cdots X_N ^{-2m_N (i)}X_M ^{-1}(X_M - X_M ^{-1})} .
\end{align*}
From Proposition~\ref{pro:sp_am} we get the following formula:
\begin{equation*}
\frac{\langle L| \mathscr{S}(n+1,n+2) \rangle }{\left | \mathscr{S}(n+1,n+2) \right |} = \frac{\langle L| \mathscr{S}(i,j) \rangle }{\left | \mathscr{S}(i,j) \right |}.
\end{equation*}
Similarly we can get the above formula if the regions $r_i$ and $r_j$ are as shown in Figure~\ref{fig:as1}. 

So the proof is complete. 
\end{pro}
\begin{Lem}\label{lem:pf2}
Put $F(L) = \langle L| \mathscr{S} \rangle / \left | \mathscr{S} \right |$. Then $F(L)$ is an invariant under the Reidemeister moves~II and III. Moreover $F(L)$ satisfies the following formulas for the Reidemeister move~I: $X_K ^{-1}F(L) = F(L') = F(L'')$ and $X_K F(L)= F(\tilde L) = F(\tilde {\tilde L})$,
where $L$, $L'$, $L''$, $\tilde L$ and $\tilde {\tilde L}$ differ only in one place as shown in Figure~\ref{fig:rm1_0}. 
\begin{figure}[h]
 \begin{center}
 \includegraphics[scale=1.15]{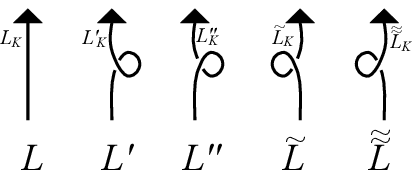}
 \caption{}
 \label{fig:rm1_0}
 \end{center}
\end{figure}
\end{Lem}
\begin{pro}
(i) Reidemeister move II:

Let $L_2$ and $L_2'$ be two link diagrams which differ only in one place as shown in Figure~\ref{fig:rm2-1}, where $L_{2,K}$ denotes the component of $L_2$ labeled as $K$. Let $\mathscr{S}$ and $\mathscr{S}'$ be the sets of states of the projections $L_2$ and $L_2'$ respectively with fixed stars as indicated in Figure~\ref{fig:rm2-1}. Then $\mathscr{S}$ has three subsets as shown in Figure~\ref{fig:rm2-1} and these subsets are denoted by $S_0$, $S_1$ and $S_2$ from left to right. On the other hand we have $|\mathscr{S}| = |\mathscr{S}'|$ for any orientations since $L_2$ and $L_2'$ have the same diagrams in the exteriors. 
\begin{figure}[h]
 \begin{center}
 \includegraphics{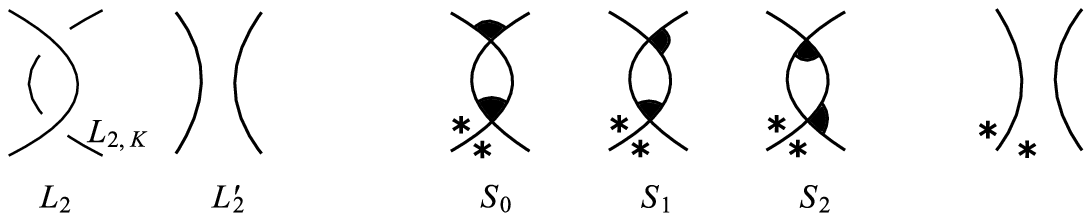}
 \caption{}
 \label{fig:rm2-1}
 \end{center}
\end{figure}
\begin{enumerate}
\item Case where the orientations are as shown in Figure~\ref{fig:rm2-2-1}: 
\\
Figure~\ref{fig:rm2-3-1} shows the variables and the signs around the vertices of $L_2$. Then we have $\frac{\langle L_2| S_1 \rangle}{-X_K^{-2}} = \frac{\langle L_2| S_2 \rangle}{X_K ^{-2}}$ since $S_1$ and $S_2$ have the same blank regions. Here a blank region means a region which has the marker outside the picture. In other words a blank region is a region which does not have markers or $\ast$ in the figure. 
Moreover we have $\langle L_2| S_0 \rangle = \langle L_2'| \mathscr{S}' \rangle$ since $S_0$ and $\mathscr{S}'$ have the same blank regions. Therefore we have $\langle L_2| \mathscr{S} \rangle = \langle L_2| S_0 \rangle + \langle L_2| S_1 \rangle + \langle L_2| S_2 \rangle = \langle L_2'| \mathscr{S} \rangle$. So we have $F(L_2) = F(L_2')$. 

If the strand over $L_{2,K}$ has the reversed orientation, $L_{2}$ has the same variables except the signs. In similar way we can get $F(L_2) = F(L_2')$ noting the signs. 
\begin{figure}[h]
\begin{minipage}{0.4\hsize}
 \begin{center}
 \includegraphics{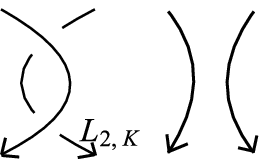}
 \caption{}
 \label{fig:rm2-2-1}
 \end{center}
 \end{minipage}
 \begin{minipage}{0.4\hsize}
 \begin{center}
 \includegraphics{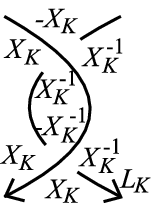}
 \caption{}
 \label{fig:rm2-3-1}
 \end{center}
\end{minipage}
\end{figure}
\item Case where the orientations are as shown in Figure~\ref{fig:rm2-2-2}:
\\
Figure~\ref{fig:rm2-3-2} shows the variables and the signs around the vertices of $L_2$. As in (i) we can prove $F(L_2) = F(L_2')$. 
\begin{figure}[h]
\begin{minipage}{0.4\hsize}
 \begin{center}
 \includegraphics{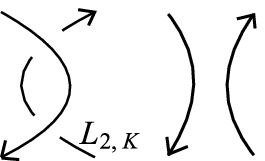}
 \caption{}
 \label{fig:rm2-2-2}
 \end{center}
 \end{minipage}
 \begin{minipage}{0.4\hsize}
 \begin{center}
 \includegraphics{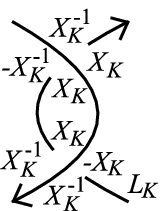}
 \caption{}
 \label{fig:rm2-3-2}
 \end{center}
\end{minipage}
\end{figure}
\end{enumerate}

(ii) Reidemeister move III:
\begin{figure}[h]
\begin{minipage}{0.4\hsize}
 \begin{center}
 \includegraphics{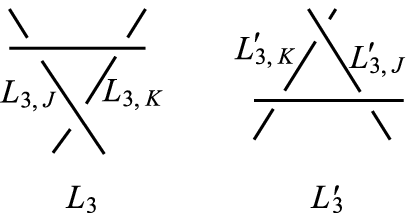}
 \caption{}
 \label{fig:rm3-1}
 \end{center}
 \end{minipage}
 \begin{minipage}{0.4\hsize}
 \begin{center}
 \includegraphics{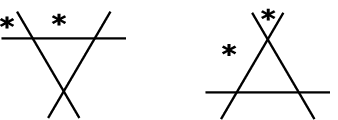}
 \caption{}
 \label{fig:rm3-2}
 \end{center}
\end{minipage}
\end{figure}

Let $L_3$ and $L_3'$ be two link diagrams which differ only in one place as shown in Figure~\ref{fig:rm3-1}. 
Let $\mathscr{S}$ and $\mathscr{S}'$ be the sets of states of the projections $L_3$ and $L_3'$ respectively with fixed stars as indicated in Figure~\ref{fig:rm3-2}. Then we have $|\mathscr{S}| = |\mathscr{S}'|$ for any orientations. On the other hand $\mathscr{S}$ and $\mathscr{S}'$ are in one-to-one correspondence for the blank regions as shown in Figure~\ref{fig:rm3-3}. Then we can prove $F(L_3)=F(L_3')$ by using the above technique for any orientation. 
\begin{figure}[H]
 \begin{center}
 \includegraphics[scale=0.9]{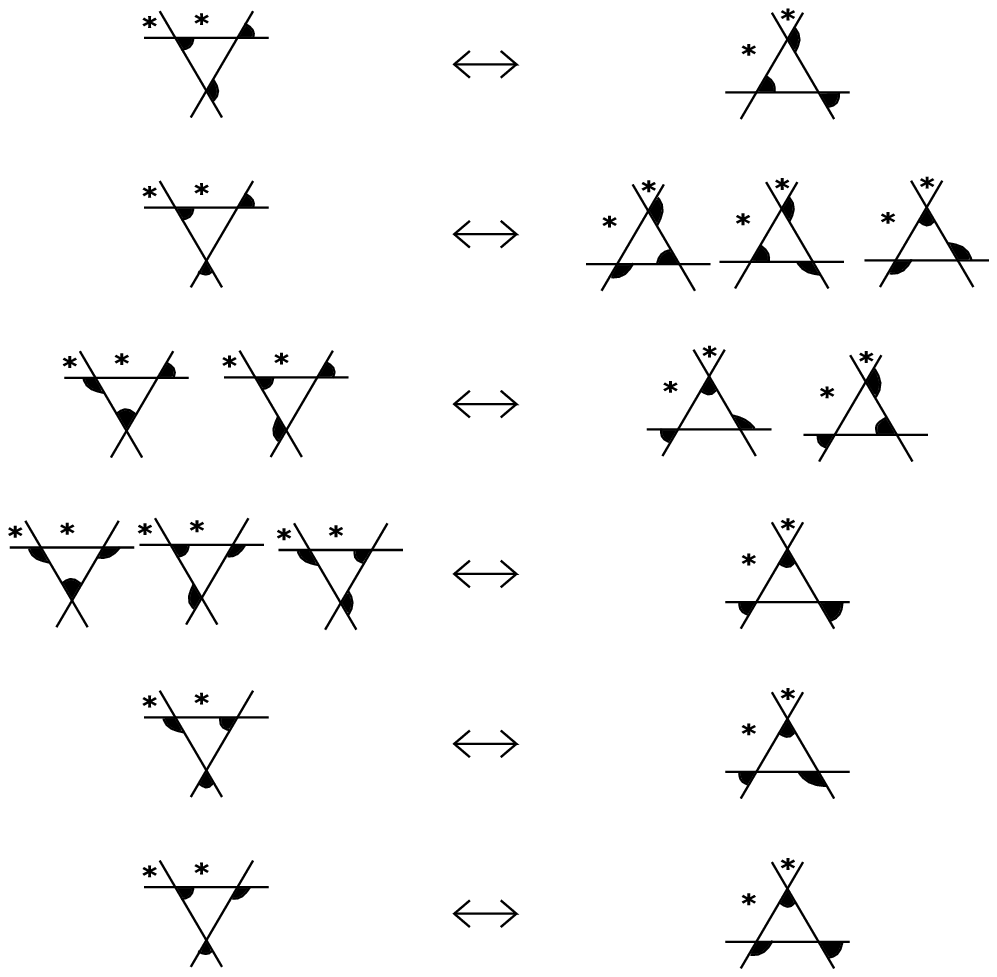}
 \caption{}
 \label{fig:rm3-3}
 \end{center}
\end{figure}

(iii) Reidemeister move I:

Let $L$, $L'$, $L''$, $\tilde L$ and $\tilde {\tilde L}$ be link diagrams which differ only in one place as shown in Figure~\ref{fig:rm1_0}. Let $\mathscr{S}$, $\mathscr{S}'$ and $\tilde {\mathscr{S}} $ be the sets of states with fixed stars as indicated in Figure~\ref{fig:rm1-1} of the projections $L$, $L'$ (or $L''$) and $\tilde L$ (or $\tilde {\tilde L}$). Then the markers are determined uniquely as indicated in Figure~\ref{fig:rm1-1}. So we have $|\mathscr{S}| = |\mathscr{S}'|$. Moreover the marker of $\mathscr{S}'$ indicates $X_K ^{-1}$ for both of $L'$ and $L''$. Therefore we have $X_K ^{-1}F(L) = F(L') = F(L'')$. On the other hand the marker of $\tilde {\mathscr{S}}$ indicates $X_K $ for both of $\tilde L$ and $\tilde {\tilde L}$. Therefore we have $X_K F(L) = F(L') = F(L'')$. 
\begin{figure}[h]
 \begin{center}
 \includegraphics[scale=1.25]{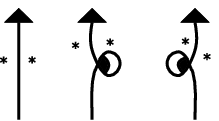}
 \caption{}
 \label{fig:rm1-1}
 \end{center}
\end{figure}

This completes the proof of Lemma~\ref{lem:pf2}. 
\end{pro}
The potential function is an invariant under the Reidemeister moves~II and III from Lemmas~\ref{lem:pf1} and \ref{lem:pf2}, and that the curvature is an invariant under these moves. Moreover the potential function is an invariant under the Reidemeister move~I from Lemmas~\ref{lem:pf2} and that the curvature counts how many times the sublink rotates counter-clockwise. So the proof of Theorem~\ref{pro:pf} is now complete properties. 
%
%
\section{Axioms for the Conway potential function}
In Section~3 we defined the potential function for links with labeled strands. 
For a finite set $\Lambda $ let $\mu :\{ 1,2,\dots ,N \} \rightarrow \Lambda $ be a surjection. A colored link $L=L_1 \cup L_2 \cup \dots \cup L_N$ with colors in $\Lambda $ is a link where each component $L_K$ is given color $\mu (K)$. Putting $\mu _K = \mu (K)$, the potential function obtained from $\square _L$ by replacing $X_K$ with $X_{\mu _K}$ is an invariant for a colored link. In this section we will show some formulas for Kauffman's potential function for colored links, and show that Kauffman's potential function equals Hartley's potential function. 

\begin{figure}[h]
 \begin{center}
 \includegraphics[scale=0.712]{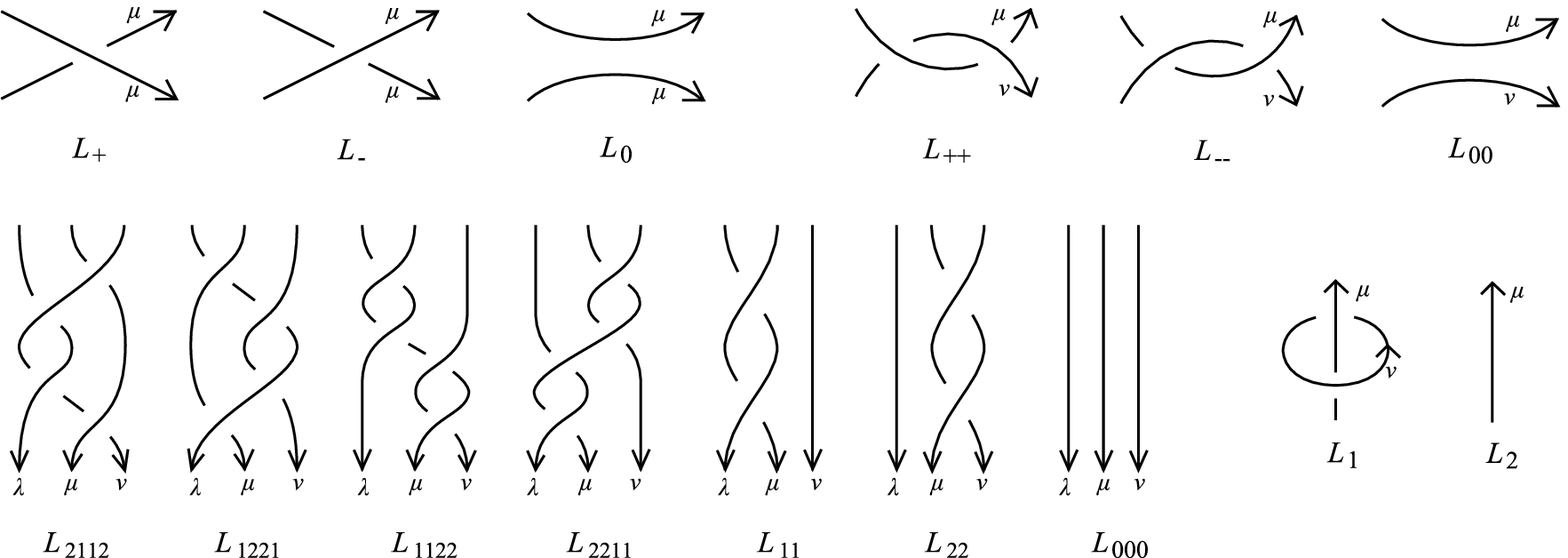}
 \caption{}
 \label{fig:pfam}
 \end{center}
\end{figure}

Let $\nabla _L$ be the potential function defined by Hartley in \cite{MR727708}. J.~Murakami showed that $\nabla _L$ can be calculated by using the following six axioms, where letters $\lambda$, $\mu$ and $\nu$ denote colors of the components. 
\begin{enumerate}
\item Let $L_{+}$, $L_{-}$ and $L_{0}$ be three links which differ only in one place as shown in Figure~\ref{fig:pfam}. Then the potential function $\nabla $ satisfies 
\begin{center}
$\nabla _{L_+} - \nabla_{L_-} = (X_\mu - X_\mu^{-1})\nabla_{L_0} $. 
\end{center}
\item Let $L_{++}$, $L_{--}$ and $L_{00}$ be three links which differ only in one place as shown in Figure~\ref{fig:pfam}. Then the potential function $\nabla $ satisfies 
\begin{center}
$\nabla _{L_{++}} - \nabla_{L_{--}} = (X_\mu X_\nu + X_\mu ^{-1}X_\nu ^{-1})\nabla_{L_{00}} $. 
\end{center}
\item Let $L_{2112}$, $L_{1221}$, $L_{1122}$, $L_{2211}$, $L_{11}$, $L_{22}$ and $L_{000}$ be seven links which differ only in one place as shown in Figure~\ref{fig:pfam}. Putting $g_+(x)=x+x^{-1}$ and $g_-(x)=x-x^{-1}$, the potential function $\nabla $ satisfies 
\begin{multline} \label{eq:pfam_3}
  g_+(X_\lambda )g_-(X_\mu )\nabla _{L_{2112}} - g_-(X_\mu )g_+(X_\nu )\nabla _{L_{1221}} \\
- g_-(X_\lambda ^{-1}X_\nu )(\nabla _{L_{1122}}+ \nabla _{L_{2211}}) + g_-(X_\lambda ^{-1}X_\mu X_\nu )g_+(X_\nu )\nabla _{L_{11}} \\
- g_+(X_\lambda )g_-(X_\lambda X_\mu X_\nu ^{-1})\nabla _{L_{22}} - g_-(X_\lambda ^{-2}X_\nu ^2)\nabla _{L_{000}} = 0. 
\end{multline}
\item For a trivial knot $L$ with color $\mu$, $\nabla_{L} = \frac{1}{X_\mu - X_\mu ^{-1}}$. 
\item Let $L_1$ and $L_2$ be two links which differ only in one place as shown in Figure~\ref{fig:pfam}. Then the potential function $\nabla $ satisfies 
\begin{center}
$\nabla _{L_1} = (X_\mu - X_\mu ^{-1})\nabla_{L_2} $. 
\end{center}
\item For a split link $L$, $\nabla_{L} = 0$. 
\end{enumerate}
\begin{Thm}[\cite{MR1177442,MR1197048}]\label{thm:jax}
The above axioms (i)--(vi) determine Hartley's potential function. 
\end{Thm}

We can show that Kauffman's potential function equals Hartley's potential function by using the above axioms. 
\begin{Thm}\label{thm:kh0}
Kauffman's potential function is equal to Hartley's potential function. 
\end{Thm}
\begin{pro}
From Theorem~\ref{thm:jax} Kauffman's potential function equals Hartley's potential function if $\square $ satisfies the above six axioms. In \cite{MR712133} Kauffman showed the axioms except (iii). So we discuss the axiom~(iii). 

Let $\mathscr{S}_{2112}$, $\mathscr{S}_{1221}$, $\mathscr{S}_{1122}$, $\mathscr{S}_{2211}$, $\mathscr{S}_{11}$, $\mathscr{S}_{22}$ and $\mathscr{S}_{000}$ be the sets of states with fixed stars as indicated in Figure~\ref{fig:pfam_3s_0} of the projections $L_{2112}$, $L_{1221}$, $L_{1122}$, $L_{2211}$, $L_{11}$, $L_{22}$ and $L_{000}$ respectively. 
\begin{figure}[h]
 \begin{center}
 \includegraphics[scale=0.75]{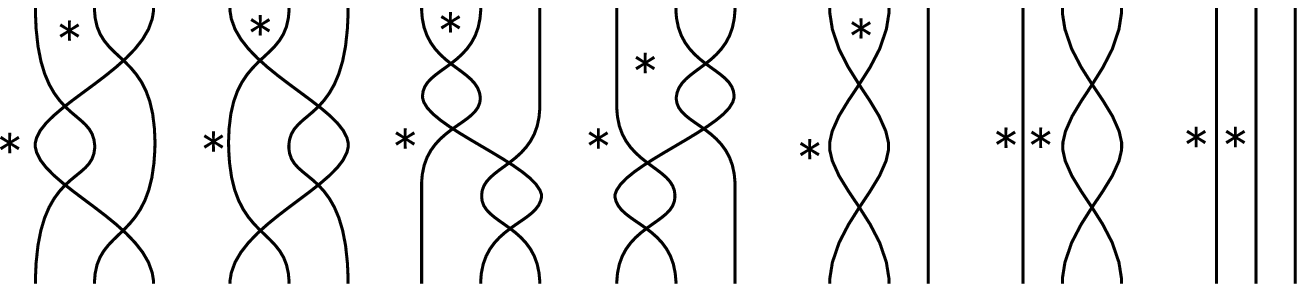}
 \caption{}
 \label{fig:pfam_3s_0}
 \end{center}
\end{figure}
Then we have $|\mathscr{S}_{2112}|=|\mathscr{S}_{1221}|=|\mathscr{S}_{1122}|=|\mathscr{S}_{2211}|=|\mathscr{S}_{11}|=|\mathscr{S}_{22}|=|\mathscr{S}_{000}|$. Moreover all the links have the same curvatures. So we will show that $\langle L| \mathscr{S} \rangle$ satisfies the axiom~(iii). 

Now we can see that the number of the regions without stars is two more than the number of the vertices for each projections in Figure~\ref{fig:pfam_3s_0}. In other words each states has two blank regions. We can classify $\mathscr{S}_{2112}$, $\mathscr{S}_{1221}$, $\mathscr{S}_{1122}$, $\mathscr{S}_{2211}$, $\mathscr{S}_{11}$, $\mathscr{S}_{22}$ and $\mathscr{S}_{000}$ as shown in Figure~\ref{fig:pfam_3s_1}--\ref{fig:pfam_3s_6}, where $\circ $ means a blank region and a dotted line expediently divides the region into a region with a marker outside the picture and a region without markers. 
\begin{figure}[H]
 \begin{center}
 \includegraphics[scale=0.75]{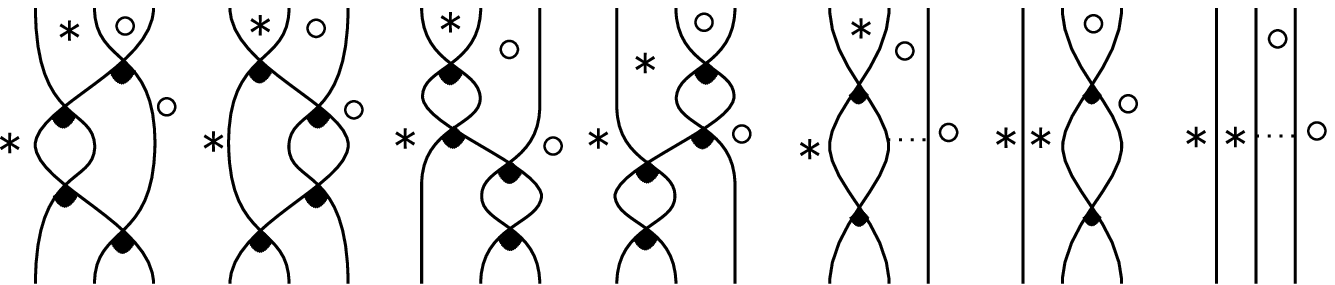}
 \caption{}
 \label{fig:pfam_3s_1}
 \end{center}
\end{figure}
\begin{figure}[H]
 \begin{center}
 \includegraphics[scale=0.75]{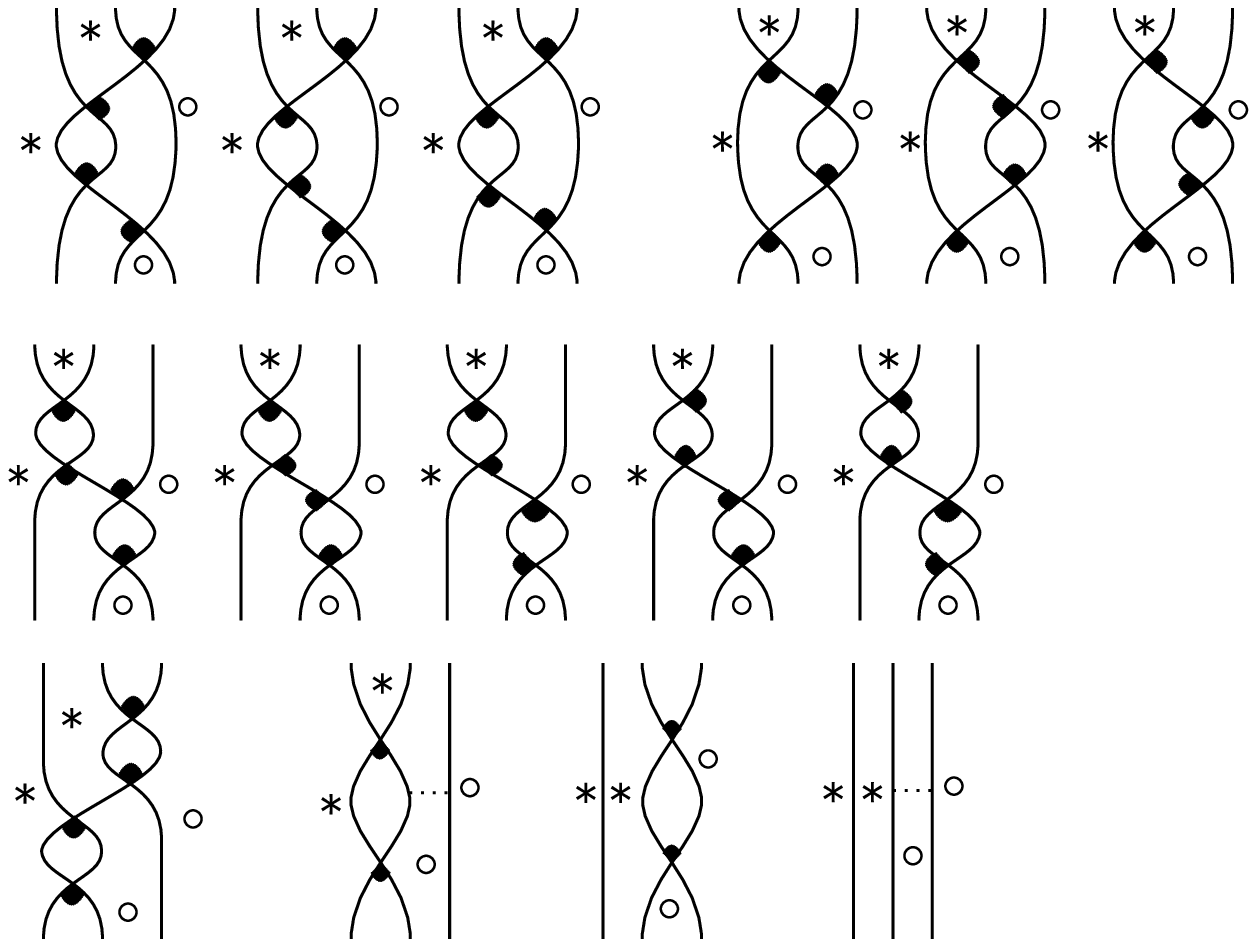}
 \caption{}
 \label{fig:pfam_3s_2}
 \end{center}
\end{figure}
\begin{figure}[H]
 \begin{center}
 \includegraphics[scale=0.75]{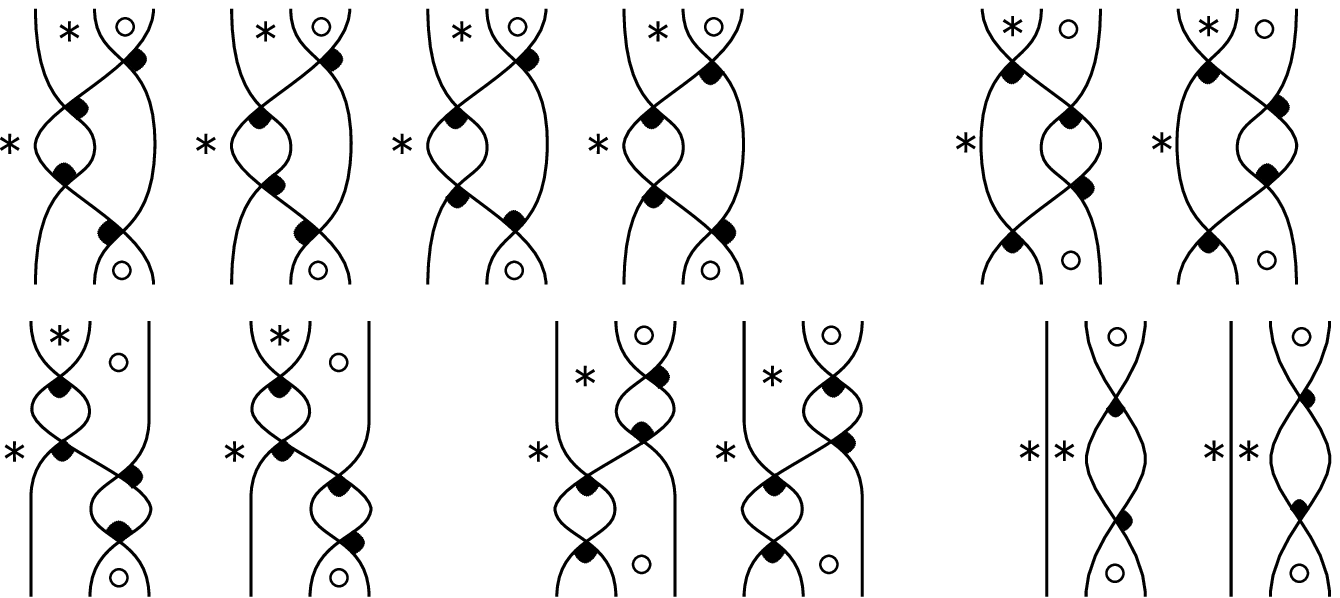}
 \caption{}
 \label{fig:pfam_3s_3}
 \end{center}
\end{figure}
\begin{figure}[H]
 \begin{center}
 \includegraphics[scale=0.75]{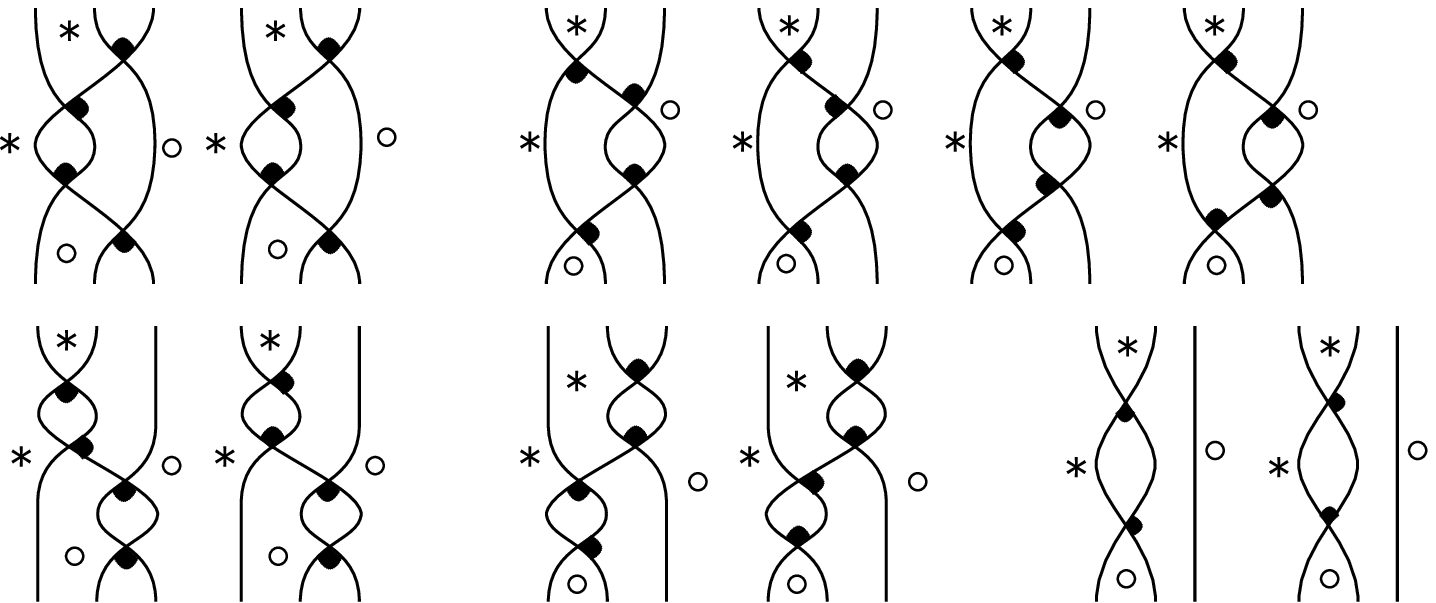}
 \caption{}
 \label{fig:pfam_3s_4}
 \end{center}
\end{figure}
\begin{figure}[H]
 \begin{center}
 \includegraphics[scale=0.75]{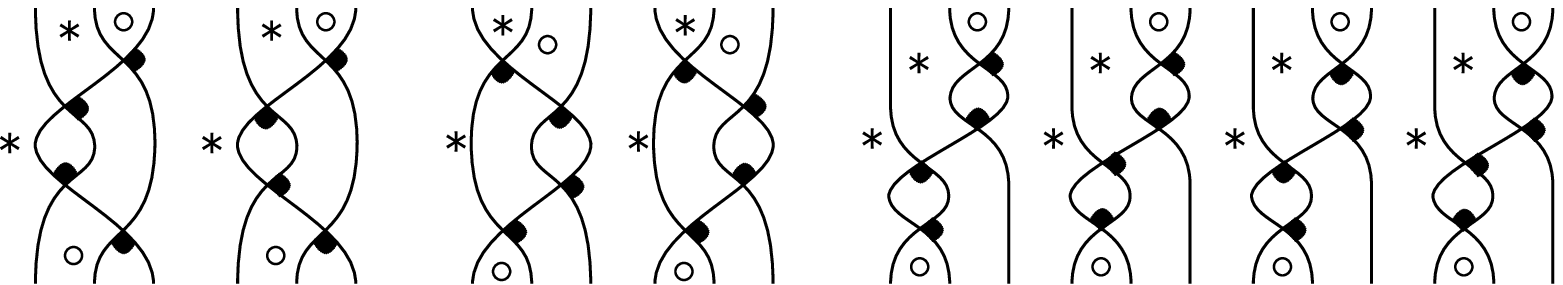}
 \caption{}
 \label{fig:pfam_3s_5}
 \end{center}
\end{figure}
\begin{figure}[H]
 \begin{center}
 \includegraphics[scale=0.75]{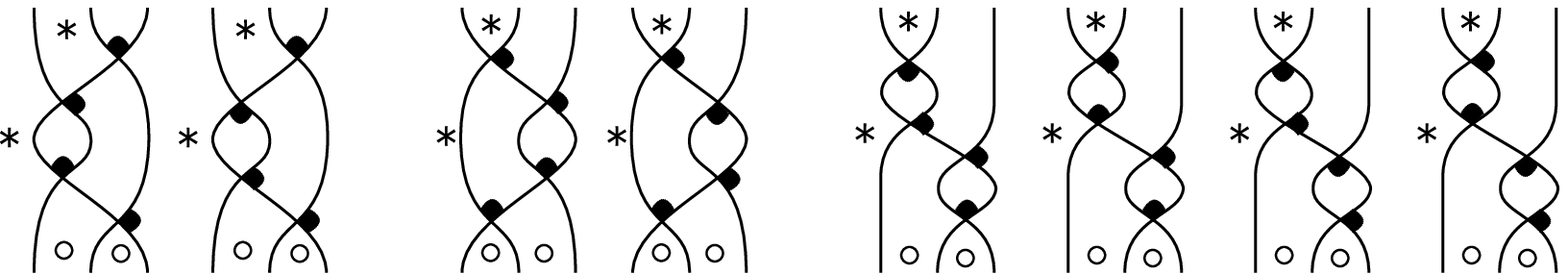}
 \caption{}
 \label{fig:pfam_3s_6}
 \end{center}
\end{figure}
Let $\mathscr{S}_{2112}^{1}$, $\mathscr{S}_{1221}^{1}$, $\mathscr{S}_{1122}^{1}$, $\mathscr{S}_{2211}^{1}$, $\mathscr{S}_{11}^{1}$, $\mathscr{S}_{22}^{1}$ and $\mathscr{S}_{000}^{1}$ be the sets of states corresponding to Figure~\ref{fig:pfam_3s_1}. Let $\mathscr{S}_{2112}^{2}$, $\mathscr{S}_{1221}^{2}$, $\mathscr{S}_{1122}^{2}$, $\mathscr{S}_{2211}^{2}$, $\mathscr{S}_{11}^{2}$, $\mathscr{S}_{22}^{2}$ and $\mathscr{S}_{000}^{2}$ be the sets of states corresponding to Figure~\ref{fig:pfam_3s_2}. Let $\mathscr{S}_{2112}^{3}$, $\mathscr{S}_{1221}^{3}$, $\mathscr{S}_{1122}^{3}$, $\mathscr{S}_{2211}^{3}$ and $\mathscr{S}_{22}^{3}$ be the sets of states corresponding to Figure~\ref{fig:pfam_3s_3}, and $\mathscr{S}_{11}^{3}$ and $\mathscr{S}_{000}^{3}$ be empty sets. Let $\mathscr{S}_{2112}^{4}$, $\mathscr{S}_{1221}^{4}$, $\mathscr{S}_{1122}^{4}$, $\mathscr{S}_{2211}^{4}$ and $\mathscr{S}_{11}^{4}$ be the sets of states corresponding to Figure~\ref{fig:pfam_3s_4}, and $\mathscr{S}_{22}^{4}$ and $\mathscr{S}_{000}^{4}$ be empty sets. Let $\mathscr{S}_{2112}^{5}$, $\mathscr{S}_{1221}^{5}$ and $\mathscr{S}_{2211}^{5}$ be the sets of states corresponding to Figure~\ref{fig:pfam_3s_5}, and $\mathscr{S}_{1122}^{5}$, $\mathscr{S}_{11}^{5}$, $\mathscr{S}_{22}^{5}$ and $\mathscr{S}_{000}^{5}$ be empty sets. Let $\mathscr{S}_{2112}^{6}$, $\mathscr{S}_{1221}^{6}$ and $\mathscr{S}_{1122}^{6}$ be the sets of states corresponding to Figure~\ref{fig:pfam_3s_6}, and $\mathscr{S}_{11}^{3}$ and $\mathscr{S}_{2211}^{6}$, $\mathscr{S}_{11}^{6}$, $\mathscr{S}_{22}^{6}$ and $\mathscr{S}_{000}^{6}$ be empty sets. 
\begin{enumerate}
\item For $\mathscr{S}_{2112}^{1}$, $\mathscr{S}_{1221}^{1}$, $\mathscr{S}_{1122}^{1}$, $\mathscr{S}_{2211}^{1}$, $\mathscr{S}_{11}^{1}$, $\mathscr{S}_{22}^{1}$ and $\mathscr{S}_{000}^{1}$:
\\
Figure~\ref{fig:pfam_3v} shows the variables and the signs around the vertices of $L_{2112}$, $L_{1221}$, $L_{1122}$, $L_{2211}$, $L_{11}$, $L_{22}$ and $L_{000}$. 
\begin{figure}[H]
 \begin{center}
 \includegraphics[scale=0.95]{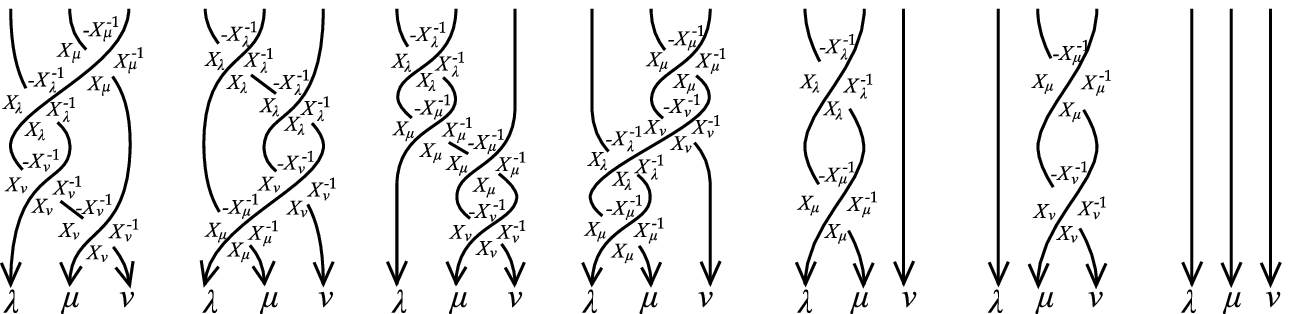}
 \caption{}
 \label{fig:pfam_3v}
 \end{center}
\end{figure}
Now the states of $\mathscr{S}_{\square}^{1}$ have the same blank regions. So we have the following formula:
\begin{equation*}
 \begin{split}
&\frac{\langle L_{2112}| \mathscr{S}_{2112}^{1} \rangle}{X_\lambda X_\mu X_\nu ^{2}} =\frac{\langle L_{1221}| \mathscr{S}_{1221}^{1} \rangle}{X_\lambda ^{2}X_\mu X_\nu } = \frac{\langle L_{1122}| \mathscr{S}_{1122}^{1} \rangle}{X_\lambda X_\mu ^{2}X_\nu } =\frac{\langle L_{2211}| \mathscr{S}_{2211}^{1} \rangle}{X_\lambda X_\mu ^{2}X_\nu } \\
=&\frac{\langle L_{11}| \mathscr{S}_{11}^{1} \rangle}{X_\lambda X_\mu } = \frac{\langle L_{22}| \mathscr{S}_{22}^{1} \rangle}{X_\mu X_\nu } =\langle L_{000}| \mathscr{S}_{000}^{1} \rangle.
 \end{split}
\end{equation*}
Therefore when we replace $\nabla _{L_\square}$ with $\langle L_\square | \mathscr{S}_{\square}^{1} \rangle$ in the left-hand side of \eqref{eq:pfam_3}, we can calculate the value as follows:
\begin{equation*}
  \begin{split}
  &g_+(X_\lambda )g_-(X_\mu )\langle L_{2112}| \mathscr{S}_{2112}^{1} \rangle - g_-(X_\mu )g_+(X_\nu )\langle L_{1221}| \mathscr{S}_{1221}^{1} \rangle \\
  &-g_-(X_\lambda ^{-1}X_\nu )(\langle L_{1122}| \mathscr{S}_{1122}^{1} \rangle + \langle L_{2211}| \mathscr{S}_{2211}^{1} \rangle ) \\ 
  &+ g_-(X_\lambda ^{-1}X_\mu X_\nu )g_+(X_\nu )\langle L_{11}| \mathscr{S}_{11}^{1} \rangle - g_+(X_\lambda )g_-(X_\lambda X_\mu X_\nu ^{-1})\langle L_{22}| \mathscr{S}_{22}^{1}\rangle \\ 
  &- g_-(X_\lambda ^{-2}X_\nu ^2)\langle L_{000}| \mathscr{S}_{000}^{1} \rangle \\
  =& \bigl{[ }g_+(X_\lambda )g_-(X_\mu )X_\lambda X_\mu X_\nu ^{2}- g_-(X_\mu )g_+(X_\nu )X_\lambda ^{2}X_\mu X_\nu \\
  &- g_-(X_\lambda ^{-1}X_\nu )(X_\lambda X_\mu ^{2}X_\nu + X_\lambda X_\mu ^{2}X_\nu ) + g_-(X_\lambda ^{-1}X_\mu X_\nu )g_+(X_\nu )X_\lambda X_\mu \\
  &- g_+(X_\lambda )g_-(X_\lambda X_\mu X_\nu ^{-1})X_\mu X_\nu - g_-(X_\lambda ^{-2}X_\nu ^2) \bigr{] }\langle L_{000}| \mathscr{S}_{000}^{1} \rangle \\
  =&\bigl{[ } X_\lambda X_\mu X_\nu \Bigl( g_-(X_\mu )\bigl( g_+(X_\lambda )X_\nu - g_+(X_\nu )X_\lambda \bigr) - 2g_-(X_\lambda ^{-1}X_\nu )X_\mu \Bigr) \\
  &+ X_\mu \bigl( g_-(X_\lambda ^{-1} X_\mu X_\nu )g_+(X_\nu )X_\lambda - g_+(X_\lambda )g_-(X_\lambda X_\mu X_\nu ^{-1})X_\nu \bigr) \\
& - g_-(X_\lambda ^{-2}X_\nu ^2)\bigr{] } \langle L_{000}| \mathscr{S}_{000}^{1} \rangle .
  \end{split}
\end{equation*}
Since the first term in the square bracket is
\begin{equation*}
  \begin{split}
  &X_\lambda X_\mu X_\nu \Bigl( g_-(X_\mu )\bigl( X_\lambda ^{-1}X_\nu - X_\lambda X_\nu ^{-1}\bigr) - 2g_-(X_\lambda ^{-1}X_\nu )X_\mu \Bigr) \\
  =& X_\lambda X_\mu X_\nu \Bigl( \bigl( X_\mu - X_\mu ^{-1}\bigr) g_-(X_\lambda ^{-1}X_\nu) - 2g_-(X_\lambda ^{-1}X_\nu )X_\mu \Bigr) \\
  =& g_-(X_\lambda ^{-1}X_\nu ) X_\lambda X_\mu X_\nu \bigl( -X_\mu ^{-1} -X_\mu \bigr) \\
  =& -g_-(X_\lambda ^{-1}X_\nu )g_+(X_\mu )X_\lambda X_\mu X_\nu ,
  \end{split}
\end{equation*}
and the second term is
\begin{equation*}
  \begin{split}
  & X_\mu \bigl( (X_\lambda ^{-1}X_\mu X_\nu - X_\lambda X_\mu ^{-1}X_\nu ^{-1})(X_\lambda X_\nu + X_\lambda X_\nu ^{-1}) \\
  &- (X_\lambda X_\mu X_\nu ^{-1}- X_\lambda ^{-1}X_\mu ^{-1}X_\nu )(X_\lambda X_\nu + X_\lambda ^{-1}X_\nu ) \bigr) \\
  =& X_\mu \bigl( X_\mu X_\nu ^2 -X_\lambda ^2X_\mu ^{-1} +X_\mu -X_\lambda ^2X_\mu ^{-1}X_\nu ^{-2} \\
  &-X_\lambda ^2X_\mu + X_\mu ^{-1}X_\nu ^2 -X_\mu +X_\lambda ^{-2}X_\mu ^{-1}X_\nu ^2 \bigr) \\
  =& X_\mu ^{2}X_\nu ^{2} -X_\lambda ^{-2} -X_\lambda ^{2}X_\nu ^{-2} -X_\lambda ^{2}X_\mu ^{2} +X_\nu ^{2} +X_\lambda ^{-2}X_\nu ^{2} \\
  =& X_\lambda X_\mu X_\nu \bigl(X_\lambda ^{-1}X_\mu X_\nu -X_\lambda X_\mu ^{-1}X_\nu ^{-1} -X_\lambda X_\mu X_\nu ^{-1} +X_\lambda ^{-1}X_\mu ^{-1}X_\nu \bigr) \\
  &-X_\lambda ^{2}X_\nu ^{-2} +X_\lambda ^{-2}X_\nu ^{2} \\
  =& X_\lambda X_\mu X_\nu \bigl( g_-(X_\lambda ^{-1}X_\nu )X_\mu + g_-(X_\lambda ^{-1}X_\nu )X_\mu ^{-1}\bigr) + g_-(X_\lambda ^{-2}X_\nu ^2) \\
  =& g_-(X_\lambda ^{-1} X_\nu ) g_+(X_\mu )X_\lambda X_\mu X_\nu + g_-(X_\lambda ^{-2}X_\nu ^2) ,
  \end{split}
\end{equation*}
\eqref{eq:pfam_3} vanishes.
%
\item For $\mathscr{S}_{2112}^{2}$, $\mathscr{S}_{1221}^{2}$, $\mathscr{S}_{1122}^{2}$, $\mathscr{S}_{2211}^{2}$, $\mathscr{S}_{11}^{2}$, $\mathscr{S}_{22}^{2}$ and $\mathscr{S}_{000}^{2}$:
\\
From Figure~\ref{fig:pfam_3s_2} and \ref{fig:pfam_3v} we have the following formula:
\begin{equation*}
 \begin{split}
&\frac{\langle L_{2112}| \mathscr{S}_{2112}^{2} \rangle}{X_\lambda ^{-1} X_\mu ^{-1}} =\frac{\langle L_{1221}| \mathscr{S}_{1221}^{2} \rangle}{X_\mu X_\nu } = \frac{\langle L_{1122}| \mathscr{S}_{1122}^{2} \rangle}{X_\lambda ^{-1}X_\nu ^{-1}-X_\lambda ^{-1}X_\nu +X_\lambda X_\nu} \\
=& \frac{\langle L_{2211}| \mathscr{S}_{2211}^{2} \rangle}{X_\lambda X_\nu ^{-1}} = \frac{\langle L_{11}| \mathscr{S}_{11}^{2} \rangle}{X_\lambda X_\mu } = \frac{\langle L_{22}| \mathscr{S}_{22}^{2} \rangle}{X_\mu ^{-1}X_\nu ^{-1}} =\langle L_{000}| \mathscr{S}_{000}^{2} \rangle .
 \end{split}
\end{equation*}
Therefore when we replace $\nabla _{L_\square}$ with $\langle L_\square | \mathscr{S}_{\square}^{2} \rangle$ in the left-hand side of \eqref{eq:pfam_3}, we can show that the value equals zero. 
%
\item For $\mathscr{S}_{2112}^{3}$, $\mathscr{S}_{1221}^{3}$, $\mathscr{S}_{1122}^{3}$, $\mathscr{S}_{2211}^{3}$, $\mathscr{S}_{11}^{3}$, $\mathscr{S}_{22}^{3}$ and $\mathscr{S}_{000}^{3}$:
\\
From Figure~\ref{fig:pfam_3s_3} and \ref{fig:pfam_3v} we have the following formulas:
\begin{equation*}
  \begin{split} 
&\frac{\langle L_{2112}| \mathscr{S}_{2112}^{3} \rangle}{-X_\lambda ^{-1} X_\mu ^{-1}+X_\lambda X_\mu } =\frac{\langle L_{1221}| \mathscr{S}_{1221}^{3} \rangle}{-X_\mu X_\nu ^{-1}+X_\lambda ^{2}X_\mu X_\nu ^{-1} }= \frac{\langle L_{1122}| \mathscr{S}_{1122}^{3} \rangle}{-X_\lambda X_\nu ^{-1} + X_\lambda X_\mu ^{2} X_\nu ^{-1} } \\
  =& \frac{\langle L_{2211}| \mathscr{S}_{2211}^{3} \rangle}{-X_\lambda X_\nu ^{-1} + X_\lambda X_\mu ^{2} X_\nu ^{-1}} = \frac{\langle L_{22}| \mathscr{S}_{22}^{3} \rangle}{-X_\mu ^{-1}X_\nu ^{-1} + X_\mu X_\nu ^{-1}} ,
  \end{split}
\end{equation*}
and
\begin{equation*}
  \langle L_{11}| \mathscr{S}_{11}^{3} \rangle = \langle L_{000}| \mathscr{S}_{000}^{3} \rangle  = 0.
\end{equation*}
Therefore when we replace $\nabla _{L_\square}$ with $\langle L_\square | \mathscr{S}_{\square}^{3} \rangle$ in the left-hand side of \eqref{eq:pfam_3}, we can show that the value equals zero. 
%
\item For $\mathscr{S}_{2112}^{4}$, $\mathscr{S}_{1221}^{4}$, $\mathscr{S}_{1122}^{4}$, $\mathscr{S}_{2211}^{4}$, $\mathscr{S}_{11}^{4}$, $\mathscr{S}_{22}^{4}$ and $\mathscr{S}_{000}^{4}$:
\\
From Figure~\ref{fig:pfam_3s_4} and \ref{fig:pfam_3v} we have the following formulas:
\begin{equation*}
  \begin{split}
&\frac{\langle L_{2112}| \mathscr{S}_{2112}^{4} \rangle}{X_\lambda ^{-1} X_\mu ^{-1} - X_\lambda X_\mu ^{-1}} = \frac{\langle L_{1122}| \mathscr{S}_{1122}^{4} \rangle}{-X_\lambda ^{-1}X_\nu + X_\lambda X_\nu } \\
  =& \frac{\langle L_{2211}| \mathscr{S}_{2211}^{4} \rangle}{-X_\lambda ^{-1} X_\mu ^{-2} X_\nu ^{-1} + X_\lambda X_\mu ^{-2} X_\nu ^{-1}} = \frac{\langle L_{11}| \mathscr{S}_{11}^{3} \rangle}{-X_\lambda ^{-1}X_\mu ^{-1} + X_\lambda X_\nu ^{-1}} ,
  \end{split}
\end{equation*}
and
\begin{equation*}
  \langle L_{1221}| \mathscr{S}_{1221}^{4} \rangle = \langle L_{22}| \mathscr{S}_{22}^{4} \rangle = \langle L_{000}| \mathscr{S}_{000}^{4} \rangle  = 0.
\end{equation*}
Therefore when we replace $\nabla _{L_\square}$ with $\langle L_\square | \mathscr{S}_{\square}^{4} \rangle$ in the left-hand side of \eqref{eq:pfam_3}, we can show that the value equals zero. 
%
\item For $\mathscr{S}_{2112}^{5}$, $\mathscr{S}_{1221}^{5}$, $\mathscr{S}_{1122}^{5}$, $\mathscr{S}_{2211}^{5}$, $\mathscr{S}_{11}^{5}$, $\mathscr{S}_{22}^{5}$ and $\mathscr{S}_{000}^{5}$:
\\
From Figure~\ref{fig:pfam_3s_5} and \ref{fig:pfam_3v} we have the following formulas:
\begin{equation*}
  \begin{split}
&\frac{\langle L_{2112}| \mathscr{S}_{2112}^{5} \rangle}{-X_\lambda ^{-1} X_\mu ^{-1} + X_\lambda X_\mu ^{-1}} = \frac{\langle L_{1221}| \mathscr{S}_{1221}^{5} \rangle}{-X_\mu ^{-1}X_\nu ^{-1} + X_\lambda ^{2}X_\nu X_\mu ^{-1}X_\nu ^{-1}} \\
  =& \frac{\langle L_{2211}| \mathscr{S}_{2211}^{5} \rangle}{X_\lambda ^{-1} X_\mu ^{-2} X_\nu ^{-1} - X_\lambda ^{-1} X_\nu ^{-1} - X_\lambda X_\mu ^{-2} X_\nu ^{-1} + X_\lambda X_\nu ^{-1}} ,
  \end{split}
\end{equation*}
and
\begin{equation*}
  \langle L_{1122}| \mathscr{S}_{1122}^{5} \rangle = \langle L_{11}| \mathscr{S}_{11}^{5} \rangle = \langle L_{22}| \mathscr{S}_{22}^{5} \rangle = \langle L_{000}| \mathscr{S}_{000}^{5} \rangle  = 0.
\end{equation*}
Therefore when we replace $\nabla _{L_\square}$ with $\langle L_\square | \mathscr{S}_{\square}^{5} \rangle$ in the left-hand side of \eqref{eq:pfam_3}, we can show that the value equals zero. 
%
\item For $\mathscr{S}_{2112}^{6}$, $\mathscr{S}_{1221}^{6}$, $\mathscr{S}_{1122}^{6}$, $\mathscr{S}_{2211}^{6}$, $\mathscr{S}_{11}^{6}$, $\mathscr{S}_{22}^{6}$ and $\mathscr{S}_{000}^{6}$:
\\
From Figure~\ref{fig:pfam_3s_6} and \ref{fig:pfam_3v} we have the following formulas:
\begin{equation*}
  \begin{split}
&\frac{\langle L_{2112}| \mathscr{S}_{2112}^{6} \rangle}{X_\lambda ^{-1} X_\mu ^{-1} X_\nu ^{-2} - X_\lambda X_\mu ^{-1} X_\nu ^{-2}} = \frac{\langle L_{1221}| \mathscr{S}_{1221}^{6} \rangle}{X_\lambda ^{-2} X_\mu ^{-1}X_\nu ^{-1} - X_\mu ^{-1}X_\nu ^{-1}} \\
  =& \frac{\langle L_{1122}| \mathscr{S}_{2211}^{6} \rangle}{X_\lambda ^{-1} X_\mu ^{-2} X_\nu ^{-1} - X_\lambda ^{-1} X_\nu ^{-1} - X_\lambda X_\mu ^{-2} X_\nu ^{-1} + X_\lambda X_\nu ^{-1}} .
  \end{split}
\end{equation*}
and
\begin{equation*}
  \langle L_{2211}| \mathscr{S}_{1122}^{6} \rangle = \langle L_{11}| \mathscr{S}_{11}^{6} \rangle = \langle L_{22}| \mathscr{S}_{22}^{6} \rangle = \langle L_{000}| \mathscr{S}_{000}^{6} \rangle  = 0 ,
\end{equation*}
Therefore when we replace $\nabla _{L_\square}$ with $\langle L_\square | \mathscr{S}_{\square}^{6} \rangle$ in the left-hand side of \eqref{eq:pfam_3}, we can show that the value equals zero. 
\end{enumerate}

Therefore we showed that $\langle L| \mathscr{S} \rangle$ satisfies axiom~(iii).

Hence the proof of Theorem~\ref{thm:kh0} is now complete. 
\end{pro}
%
\bibliography{2010_7_13}
\bibliographystyle{amsplain}
\end{document}